\documentclass[12pt,twoside]{article}

\usepackage{amssymb}

\renewcommand{\leq}{\leqslant}
\renewcommand{\geq}{\geqslant}

\begin{document}
\parindent=0.3in
\parskip=0in
\baselineskip=18pt plus 1pt

\vspace*{0.1in}
\noindent
{\large \bf
Factorisation of Lie Resolvents}

\vspace{0.2in}
\noindent
R. M. BRYANT \hfill roger.bryant@manchester.ac.uk\\
{\footnotesize {\it School of Mathematics, University of Manchester,
PO Box 88, Manchester M60 1QD, UK}}

\vspace{0.2in}
\noindent
M. SCHOCKER\footnote{\ Research of second author supported by Deutsche
Forschungsgemeinschaft (DFG Scho-799).}
\hfill m.schocker@swansea.ac.uk\\
{\footnotesize {\it Department of Mathematics, University of Wales Swansea,
Singleton~Park,~Swansea~SA2~8PP,~UK}}

\vspace{0.2in}
\noindent
{\bf Abstract.} \hspace{0.1in} Let $G$ be a group, $F$ a field of prime
characteristic $p$ and $V$ a finite-dimensional $FG$-module.
Let $L(V)$ denote the free Lie algebra on $V$, regarded as
an $FG$-module, and, for each positive integer $r$, let $L^r(V)$
be the $r$th homogeneous component of $L(V)$, 
called the $r$th Lie power of $V$. In a previous paper we obtained
a decomposition of $L^r(V)$ as a direct sum of modules of the form
$L^s(W)$, where $s$ is a power of $p$. Here we derive some consequences. 
First we obtain a similar result for restricted Lie powers of $V$.
Then we consider the `Lie resolvents' $\Phi^r\,$: certain functions on the
Green ring of $FG$ which determine Lie powers up to isomorphism. 
For $k$ not divisible by $p$, we obtain the factorisation
$\Phi^{p^mk} = \Phi^{p^m} \circ \Phi^k$, separating out the key case of
$p$-power degree. Finally we study certain functions on power series
over the Green ring, denoted by 
${\bf S}^*$ and ${\bf L}^*$, 
which encode symmetric powers and Lie powers, respectively.
In characteristic $0$, ${\bf L}^*$ is the
inverse of ${\bf S}^*$. In characteristic $p$, the
composite ${\bf L}^* \circ {\bf S}^*$ maps any $p$-typical power series
to a $p$-typical power series.

\vspace{0.2in}
\noindent
{\bf Keywords:} \hspace{0.1in} free Lie algebra, Lie power,
Adams operation, Lie resolvent 

\vspace{0.1in}
\noindent
{\bf Mathematics Subject Classification (2000):}
\hspace{0.1in} 17B01, 20C07, 20C20

\vspace{0.4in}
\baselineskip=20pt plus 1pt
\noindent
{\bf 1. Introduction}
\renewcommand{\theequation}{1.\arabic{equation}}
\setcounter{equation}{0}

\noindent
Let $G$ be a group and $F$ a field. For any finite-dimensional
$FG$-module $V$ let $L(V)$ be the free Lie algebra on $V$. Then
$L(V)$ may be regarded as an $FG$-module on which each element
of $G$ acts as a Lie algebra automorphism. Furthermore, each
homogeneous component $L^r(V)$ is a finite-dimensional submodule,
called the $r$th Lie power of $V$.  We regard $L(V)$ as an
$FG$-submodule of the free associative algebra $T(V)$. Thus
$L^r(V)$ is a submodule of the $r$th homogeneous component $T^r(V)$.

When $F$ has prime characteristic $p$ we may also form the free
restricted Lie algebra $R(V)$ which again may be regarded as
an $FG$-submodule
of $T(V)$. The homogeneous component $R^r(V)$ is called the
$r$th restricted Lie power of $V$.

A general problem is to describe the modules $L^r(V)$ and $R^r(V)$
up to isomorphism. We refer to [4] and the papers cited there
for a discussion of progress on $L^r(V)$. Although some use has been
made of $R(V)$ in studying $L(V)$, results so far on the module
structure of $R^r(V)$ are rather sparse: but see [9, 10, 17].

The main result of [4], called 
the `Decomposition Theorem', reduces the study
of arbitrary Lie powers in characteristic $p$
to the study of Lie powers of
$p$-power degree. It states that, for each positive
integer $r$, there is a submodule $B_r$ of $L^r(V)$ such that
$B_r$ is a direct summand of $T^r(V)$ and, for $k$ not
divisible by $p$ and $m \geq 0$,
$$L^{p^mk}(V) = L^{p^m}(B_k) \oplus L^{p^{m-1}}(B_{pk})
\oplus \cdots \oplus L^1(B_{p^mk}).$$
In this paper we derive some consequences of the Decomposition
Theorem.

Our first result, obtained in \S2, is an analogous result
for $R(V)$. With the modules $B_r$ and $k$ and $m$ as before,
we find that
$$R^{p^mk}(V) = R^{p^m}(B_k) \oplus R^{p^{m-1}}(B_{pk})
\oplus \cdots \oplus R^1(B_{p^mk}).$$

In the remainder of the paper we turn to applications
in the Green ring $R_{FG}$.
This is the ring spanned by the isomorphism classes of
finite-dimensional $FG$-modules with addition and multiplication
coming from direct sums and tensor products, respectively.
In \S3 we describe the general framework which may be
used to study Lie powers and symmetric powers of
finite-dimensional $FG$-modules. This is mainly a summary of
material contained in [1]. In particular we
describe the relevance of the Adams operations and the
Lie resolvents. The latter are $\mathbb Z$-linear maps
$\Phi^r : R_{FG} \to R_{FG}$.
Knowledge of these maps is essentially equivalent to knowledge
of the isomorphism classes of Lie powers of finite-dimensional
$FG$-modules.  We also introduce certain functions ${\bf L}^*$
and ${\bf S}^*$
defined on formal power series in an indeterminate $t$ with
coefficients from $R_{FG}$. Here ${\bf L}^*$ encodes Lie powers
and ${\bf S}^*$ encodes symmetric powers. Furthermore,
${\bf L}^*$ and ${\bf S}^*$ have 
properties like those of the logarithm function and exponential
function, respectively.

In \S4 we obtain the main result of the paper, called the
`Factorisation Theorem'. It gives a `factorisation' of Lie
resolvents under composition: 
for every non-negative integer $m$ and every positive integer
$k$ not divisible by $p$, $\Phi^{p^mk} = \Phi^{p^m}
\circ \Phi^k$. The Lie resolvents $\Phi^k$ for $k$
not divisible by $p$ are comparatively well understood, so the
result can be interpreted as a reduction to the case
of $p$-power degree. The Factorisation Theorem was conjectured in [3]
and proved in [7] in the special case $m=1$.
It is interesting that the Witt polynomials (as used to define the
operations on Witt vectors) arise here in connection with
the Factorisation Theorem.

Finally, in \S5, we describe relations between ${\bf L}^*$
and ${\bf S}^*$.  In characteristic $p$ we
use the Factorisation Theorem to show that the composite
function ${\bf L}^* \circ {\bf S}^*$ takes any power series
of the form $Vt$ to a power series where the coefficient of
$t^r$ is $0$ unless $r$ is a power of $p$. A power series of
the latter form is called `$p$-typical'. It follows that
${\bf L}^* \circ {\bf S}^*$ maps every $p$-typical power series to
a $p$-typical power series.
In characteristic $0$ we find that ${\bf L}^*$ is
the inverse of ${\bf S}^*$. This is closely related to results of
Joyal [12] and Reutenauer [14].

\vspace{0.4in}
\noindent
{\bf 2. Decomposition of Lie powers and restricted Lie powers}
\renewcommand{\theequation}{2.\arabic{equation}}
\setcounter{equation}{0}

\noindent
Let $V$ be a finite-dimensional $FG$-module, where $F$ is a field and
$G$ is a group. We write $T(V)$ for the free associative algebra
freely generated by any $F$-basis of $V$. Thus $T(V)$ has an $F$-space
decomposition $T(V) = \bigoplus_{r \geq 0}T^r(V)$, where, for each $r$,
$T^r(V)$ is the $r$th homogeneous component of $T(V)$ and $T^1(V)$
is identified with $V$.
The action of $G$ on $V$ extends to $T(V)$ so that $G$ acts by algebra
automorphisms. Thus $T(V)$ becomes an $FG$-module, and each $T^r(V)$
is a finite-dimensional submodule.

The algebra $T(V)$ may be made into a Lie algebra by defining
$[a,b] = ab-ba$ for all $a,b \in T(V)$. If $F$ has prime characteristic
$p$ then $T(V)$ may be made into a restricted Lie algebra by taking 
the additional powering operation
$a \mapsto a^p$. The Lie subalgebra of $T(V)$ generated by $V$ is denoted
by $L(V)$ and (when $F$ has characteristic $p$) the restricted Lie
subalgebra of $T(V)$ generated by $V$ is denoted by $R(V)$. As is well
known, 
$L(V)$ is a free Lie algebra and $R(V)$ is a free restricted Lie algebra,
both freely generated by any basis of $V$ (see, for example, [13, Sections
1.2, 1.6.3 and 2.5.2]). 
For $r \geq 1$ we write $L^r(V) = T^r(V) \cap L(V)$ and
$R^r(V) = T^r(V) \cap R(V)$. Thus $L(V) = \bigoplus_{r \geq 1}L^r(V)$
and $R(V) = \bigoplus_{r \geq 1}R^r(V)$. Here $L^r(V)$ and $R^r(V)$
are $FG$-submodules of $T^r(V)$, called, respectively, the $r$th Lie
power and $r$th restricted Lie power of $V$.  

As observed in [4, \S2], if $B$ is a submodule of $T^r(V)$, 
the associative subalgebra of $T(V)$ generated by $B$ may be identified
with $T(B)$. The Lie subalgebra of $T(V)$ generated by $B$ is then the 
free Lie algebra $L(B)$ and the restricted Lie subalgebra of $T(V)$
generated by $B$ is the free restricted Lie algebra $R(B)$. We assume
such identifications in the next two theorems. The first is
[4, Theorem 4.4]. 

\vspace{0.2in}
\noindent
{\bf Theorem 2.1} (Decomposition Theorem) \hspace{0.1in}
{\it Let\/ $F$ be a field
of prime characteristic\/ $p$. Let\/ $G$ be a group and\/ $V$ a
finite-dimensional\/ $FG$-module. For each positive integer\/ $r$ there
is a submodule\/ $B_r$ of\/ $L^r(V)$ such that\/ $B_r$ is a direct
summand of\/ $T^r(V)$ and, for\/ $k$ not divisible by\/ $p$ and\/
$m \geq 0$, 
\begin{equation}
L^{p^mk}(V) = L^{p^m}(B_k) \oplus L^{p^{m-1}}(B_{pk}) \oplus
\cdots \oplus L^1(B_{p^mk}).
\end{equation}}

\vspace{0.2in}
We use this to prove a similar result for restricted Lie powers.

\vspace{0.2in}
\noindent
{\bf Theorem 2.2} \hspace{0.1in} {\it In the notation of Theorem\/} 2.1,
$$R^{p^mk}(V) = R^{p^m}(B_k) \oplus R^{p^{m-1}}(B_{pk}) \oplus
\cdots \oplus R^1(B_{p^mk}).$$

\vspace{0.2in}
\noindent
{\bf Proof:} \hspace{0.1in}
For any subset $\cal S$ of $L(V)$ and $i \geq 0$ we
write ${\cal S}^{[p^i]}$ for the set $\{ s^{p^i} : s \in {\cal S}\}$.
We shall use the fact that if $\cal X$ is any basis of $L(V)$ then
the elements $x^{p^i}$ with $x \in {\cal X}$ and $i \geq 0$ are
distinct and form a basis of $R(V)$: they are linearly independent
by the Poincar\'e--Birkhoff--Witt Theorem and it is easy to see that
they span $R(V)$. 

For each non-negative integer $i$ and each positive integer $k$ not
divisible by $p$, let ${\cal X}(i,k)$ be a basis of $L^{p^ik}(V)$. Hence
$\bigcup_{i,k}{\cal X}(i,k)$ is a basis of $L(V)$. Here, and in the
remainder of the proof, all unions are unions of disjoint sets. 
Furthermore, write
\begin{equation}
\widehat{\cal X}(i,k) = {\cal X}(0,k)^{[p^i]} \cup
{\cal X}(1,k)^{[p^{i-1}]} \cup \cdots \cup {\cal X}(i,k)^{[1]}.
\end{equation}
Considerations of degree in the basis of $R(V)$ obtained from
$\bigcup_{i,k}{\cal X}(i,k)$ show that $\widehat{\cal X}(i,k)$ is a basis
of $R^{p^ik}(V)$.

For each triple $(i,j,k)$ where $i$ and $j$ are non-negative integers
and $k$ is a positive integer not divisible by $p$, let
${\cal Y}(i,j,k)$ be any basis of $L^{p^i}(B_{p^jk})$ and write
\begin{equation}
\widehat{\cal Y}(i,j,k) = {\cal Y}(0,j,k)^{[p^i]} \cup
{\cal Y}(1,j,k)^{[p^{i-1}]} \cup \cdots \cup {\cal Y}(i,j,k)^{[1]}.
\end{equation}
Thus $\widehat{\cal Y}(i,j,k)$ is a basis of $R^{p^i}(B_{p^jk})$.

By the Decomposition Theorem, we can choose the basis ${\cal X}(i,k)$
of $L^{p^ik}(V)$ to satisfy
\begin{equation}
{\cal X}(i,k) = {\cal Y}(i,0,k) \cup {\cal Y}(i-1,1,k) \cup \cdots \cup
{\cal Y}(0,i,k).
\end{equation}
Let $m$ be a non-negative integer. It is easily verified from (2.2) and
(2.4) that $\widehat{\cal X}(m,k)$ is the union of the sets
${\cal Y}(i,j,k)^{[p^r]}$ with $i+j+r = m$. Hence, by (2.3),
$$\widehat{\cal X}(m,k) = \widehat{\cal Y}(m,0,k) \cup
\widehat{\cal Y}(m-1,1,k) \cup \cdots \cup \widehat{\cal Y}(0,m,k).$$
Since $\widehat{\cal X}(m,k)$ spans $R^{p^mk}(V)$ and
$\widehat{\cal Y}(i,m-i,k)$ spans $R^{p^i}(B_{p^{m-i}k})$, we obtain
the required result. \hfill{$\square$}

\vspace{0.4in}
\noindent
{\bf 3. Symmetric powers and Lie powers}
\renewcommand{\theequation}{3.\arabic{equation}}
\setcounter{equation}{0}

\noindent
The underlying ideas that we use are
described in [1]. However,
we shall make the treatment here as self-contained as possible. We
also formulate some of the ideas in a slightly new way,
in order to emphasise the analogies between symmetric powers and Lie
powers. 
Let $G$ be a group and $F$ a field. Let ${\rm Mod}(FG)$ be the class of
all finite-dimensional $FG$-modules and let $R_{FG}$ be the associated
Green ring (or representation ring), as defined in \S1. 

The $\mathbb Q$-algebra ${\mathbb Q} \otimes_{\mathbb Z} R_{FG}$
will be denoted by $\Gamma_{FG}$ and we regard $R_{FG}$ as a subset
of $\Gamma_{FG}$. (In [1], $\mathbb C$ was used instead
of $\mathbb Q$. But, for what is needed here, 
$\mathbb Q$ works as well as $\mathbb C$.)

Let $\Pi$ denote the power series ring over $\Gamma_{FG}$ in
an indeterminate $t$; that is, $\Pi = \Gamma_{FG}[[t]]$. 
Thus $t\Pi$ and $1 + t\Pi$ are the subsets consisting of all power series
with constant term $0$ and $1$, respectively.  For
$f \in t\Pi$ and $g \in 1 + t\Pi$ we may define
${\rm exp}(f) \in 1 + t\Pi$ and ${\rm log}(g) \in
t\Pi$ by
$${\rm exp}(f) = 1 + f + f^2/2! + f^3/3! + \cdots ,$$
$${\rm log}(g) = (g-1) - (g-1)^2/2 + (g-1)^3/3 - \cdots \: .$$
Thus we have mutually inverse functions
$${\rm exp} : t\Pi \to 1 + t\Pi,
\quad {\rm log} : 1 + t\Pi \to t\Pi.$$

If $V \in {\rm Mod}(FG)$ we also write $V$ for the element of $R_{FG}$
or $\Gamma_{FG}$ determined by $V$. Thus, for example, $T^r(V)$ as
an element of $R_{FG}$ denotes the isomorphism class of the module
$T^r(V)$ and it is equal to $V^r$.

We start with a discussion of symmetric powers and Adams operations.
Let $V \in {\rm Mod}(FG)$. We write $S(V)$ for the polynomial
algebra (free commutative associative algebra) over $F$ freely
generated by any basis of $V$. The action of $G$ on $V$ extends to
$S(V)$ so that $G$ acts by algebra automorphisms. For $r \geq 0$, the
$r$th homogeneous component $S^r(V)$ is an $FG$-submodule of $S(V)$
called the $r$th symmetric power of $V$. Here $S^1(V)$ is identified 
with $V$.

For $V \in {\rm Mod}(FG)$, let $S(V,t) \in 1 + t\Pi$
be defined by
$$
S(V,t) = 1 + S^1(V)t + S^2(V)t^2 + \cdots \: .
$$
It is well known and easy to verify that for $U,V \in {\rm Mod}(FG)$
we have
$$
S(U \oplus V,t) = S(U,t)S(V,t).
$$
It follows that there is a $\mathbb Q$-linear function $\psi :
\Gamma_{FG} \to t\Pi$ satisfying
\begin{equation}
\psi(V) = {\rm log}(S(V,t))
\end{equation}
for all $V \in {\rm Mod}(FG)$. Hence we may define a function
$${\bf S}^+ : t\Pi \longrightarrow t\Pi$$
by
\begin{equation}
{\bf S}^+(A_1t + A_2t^2 + A_3t^3 + \cdots) =
\psi(A_1) + \psi(A_2)_{t \mapsto t^2} + \psi(A_3)_{t \mapsto t^3}
+ \cdots \: ,
\end{equation}
for all $A_1,A_2,\dots \in \Gamma_{FG}$, where the subscript
$t \mapsto t^r$ denotes the operation of replacing a power series
$X_1t + X_2t^2 + \cdots\;$ by $X_1t^r + X_2t^{2r} + \cdots \:$.
The properties of ${\bf S}^+$ given in the following lemma are
readily obtained from the definition.

\vspace{0.2in}
\noindent
{\bf Lemma 3.1} \hspace{0.1in}
{\it The function\/ ${\bf S}^+$ is\/ $\mathbb Q$-linear
and\/ ${\bf S}^+(t^r\Pi) \subseteq t^r\Pi$, for all\/ $r \geq 1$.
For all\/ $f \in t\Pi$ and all\/ $r \geq 1$,
${\bf S}^+(f_{t \mapsto t^r}) = {\bf S}^+(f)_{t \mapsto t^r}$.
If\/ $f_i \in t^i\Pi$ for\/ $i=1,2,\dots$, then
${\bf S}^+(\sum f_i) = \sum {\bf S}^+(f_i)$.}

\vspace{0.2in}
We now define a function
$${\bf S}^* : t\Pi \longrightarrow 1 + t\Pi$$
as the composite
\begin{equation}
{\bf S}^* = {\rm exp} \circ {\bf S}^+.
\end{equation}
It is clear from (3.1) and (3.2) that 
${\bf S}^*(Vt) = S(V,t)$
for all $V \in {\rm Mod}(FG)$. Also, since ${\bf S}^+$ is additive, 
$$
{\bf S}^*(f + g) = {\bf S}^*(f){\bf S}^*(g)
$$
for all $f,g \in t\Pi$.

For $A \in \Gamma_{FG}$ we may define $\psi^1(A),
\psi^2(A),\dots \in \Gamma_{FG}$ by the equation
\begin{equation}
{\bf S}^+(At) = \psi(A) = \psi^1(A)t + \frac{1}{2}\psi^2(A)t^2
+ \frac{1}{3}\psi^3(A)t^3 + \cdots \: .
\end{equation}
Thus we obtain
$\mathbb Q$-linear functions $\psi^r : \Gamma_{FG} \to \Gamma_{FG}$.
These functions were denoted by $\psi_S^r$ in [1, \S5]: they are 
the {\it Adams operations} on $\Gamma_{FG}$ 
determined by symmetric powers. 
(Adams operations, different in general, can also be defined
using exterior powers.) 
By [1, \S5], if $V \in {\rm Mod}(FG)$
then $\psi^r(V) \in R_{FG}$. Thus $\psi^r$ restricts to a
$\mathbb Z$-linear function $\psi^r : R_{FG} \to R_{FG}$.

If $V \in {\rm Mod}(FG)$ then, by (3.1) and (3.4),
$$
{\rm log}(S(V,t)) = \psi^1(V)t + \frac{1}{2}\psi^2(V)t^2 + \cdots \: .
$$
This equation shows that symmetric powers may be expressed in terms
of Adams operations and vice versa.

\vspace{0.2in}
\noindent
{\bf Lemma 3.2} \hspace{0.1in}
{\it If\/ $r$ and\/ $s$ are positive integers such
that\/ $r$ is not divisible by the characteristic of\/ $F$ then\/
$\psi^r$ is an algebra endomorphism of\/ $\Gamma_{FG}$ and
$$\psi^{rs} = \psi^r \circ \psi^s.$$
In particular, if\/ $F$ has prime characteristic\/ $p$ then
$$\psi^{p^mk} = \psi^k \circ \psi^{p^m}$$
for all non-negative integers\/ $m$ and all positive integers\/ $k$
not divisible by\/ $p$.}

\vspace{0.2in}
\noindent
{\bf Proof:} \hspace{0.1in} See [1, Theorem 5.4]. \hfill{$\square$}

\vspace{0.2in}
All the results stated above for symmetric powers have analogues
for exterior powers. However, we shall not need exterior powers
for the applications in this paper. 

We now turn to Lie powers and Lie resolvents.
We shall summarise the general theory
described in [1] and
formulate some of it in a slightly different way. 
The theory is based on a function
$${\cal L}_{FG} : t\Pi \longrightarrow t\Pi$$
called the {\it Lie module function}. When $F$ and
$G$ are understood we write $\cal L$ instead of ${\cal L}_{FG}$. 
The properties of this function are described in [1, \S2 and \S3]. In
particular, by [1, Lemma 3.1], 
\begin{equation}
{\cal L}(Vt^r) = L^1(V)t^r + L^2(V)t^{2r} + \cdots \: ,
\end{equation}
for all $V \in {\rm Mod}(FG)$ and $r \geq 1$, and
\begin{equation}
{\cal L}(f) + {\cal L}(g) = {\cal L}(f + g -fg),
\end{equation}
for all $f,g \in t\Pi$.

There are advantages in considering a modification of ${\cal L}$,
namely the function
$${\bf L}^* : 1 + t\Pi \longrightarrow t\Pi$$
defined by
\begin{equation}
{\bf L}^*(f) = - {\cal L}(1-f)
\end{equation}
for all $f \in 1 + t\Pi$. It follows from (3.5), (3.6) and
(3.7) that
$$
-{\bf L}^*(1-Vt) = L^1(V)t + L^2(V)t^2 + \cdots \: ,
$$
for all $V \in {\rm Mod}(FG)$, and 
$$
{\bf L}^*(fg) = {\bf L}^*(f) + {\bf L}^*(g),
$$
for all $f,g \in 1 + t\Pi$.

We may define a function
$${\bf L}^+ : t\Pi \longrightarrow t\Pi$$
by ${\bf L}^+ = {\bf L}^* \circ {\rm exp}$. Thus ${\bf L}^+$ is
additive and
\begin{equation}
{\bf L}^* = {\bf L}^+ \circ {\rm log}.
\end{equation}

Let
${\rm Exp} : t\Pi \to t\Pi$ be defined by
${\rm Exp}(f) = 1 - {\rm exp}(-f)$ for all $f \in t\Pi$.
Since
$${\bf L}^+(f) = -{\bf L}^+(-f) = - {\bf L}^*({\rm exp}(-f))
= {\cal L}(1-{\rm exp}(-f)),$$
we obtain 
\begin{equation}
{\bf L}^+ = {\cal L} \circ {\rm Exp}.
\end{equation}
Additional properties of
${\bf L}^+$ are given in [1, \S3], where ${\bf L}^+$ is denoted by
${\cal F}_{FG}$. These properties yield the following result. 

\vspace{0.2in}
\noindent
{\bf Lemma 3.3} \hspace{0.1in}
{\it The function\/ ${\bf L}^+$ has the properties
of\/ ${\bf S}^+$ stated in Lemma\/} 3.1.

\vspace{0.2in}
\noindent 
In particular ${\bf L}^+$ is $\mathbb Q$-linear.

For $A \in \Gamma_{FG}$ we may define
$\Phi^1(A), \Phi^2(A), \dots \in \Gamma_{FG}$ by the
equation
\begin{equation}
{\bf L}^+(At) =
\Phi^1(A)t + \frac{1}{2}\Phi^2(A)t^2 + \frac{1}{3}
\Phi^3(A)t^3 + \cdots \: .
\end{equation}
Thus we obtain $\mathbb Q$-linear functions $\Phi^r : \Gamma_{FG} \to
\Gamma_{FG}$. These functions were considered in
[1], where they were written $\Phi_{FG}^r$
and called the {\it Lie resolvents\/} for $G$
over $F$. They restrict to $\mathbb Z$-linear functions
$\Phi^r : R_{FG} \to R_{FG}$.  Lie powers may be expressed
in terms of Lie resolvents and vice versa by means of the
following expressions given in [1, Corollary 3.3]:
for all $V \in {\rm Mod}(FG)$ and every positive integer $r$,
\begin{equation}
L^r(V) = \frac{1}{r} \sum_{d \mid r} \Phi^d(V^{r/d})
\end{equation}
and
\begin{equation}
\Phi^r(V) = \sum_{d \mid r} \mu(r/d) \, d \, L^d(V^{r/d}),
\end{equation}
where $\mu$ denotes the M\"obius function.

We conclude this section with two further lemmas.

\vspace{0.2in}
\noindent
{\bf Lemma 3.4} \hspace{0.1in}  {\it If\/ $r$ is not divisible
by the characteristic of\/ $F$, then
$$\Phi^r = \mu(r)\psi^r.$$}

\vspace{0.01in}
\noindent
{\bf Proof:} \hspace{0.1in} See [1, Corollary 6.2]. \hfill{$\square$}

\vspace{0.2in}
\noindent
{\bf Lemma 3.5} \hspace{0.1in}
{\it Let\/ $f \in t\Pi$, where\/ $f = A_1t +
A_2t^2 + \cdots\:$. Then 
$${\bf L}^+(f) = \sum_{r \geq 1} \biggl(
\sum_{d \mid r} \frac{1}{d} \Phi^d(A_{r/d})\biggr) t^r$$
and
$${\bf S}^+(f) = \sum_{r \geq 1} \biggl(
\sum_{d \mid r} \frac{1}{d} \psi^d(A_{r/d})\biggr) t^r.$$}

\vspace{0.2in}
\noindent
{\bf Proof:} \hspace{0.1in}
By (3.10) and (3.4), and in the terminology of [1, \S2],
the functions $\frac{1}{r}\Phi^r$
and $\frac{1}{r}\psi^r$ are the `components' of ${\bf L}^+$ and
${\bf S}^+$, respectively.  Hence
the result follows from [1, (2.5)]. \hfill{$\square$}

\newpage
\noindent
{\bf 4. The Factorisation Theorem}
\renewcommand{\theequation}{4.\arabic{equation}}
\setcounter{equation}{0}

\noindent
Let $F$ be an infinite field and $n$ a positive integer.
We write $G(n) = {\rm GL}(n,F)$.  We require some facts about polynomial
modules and their characters. Most of these have already been
collected in [1, \S5], so we use this as a convenient reference:
but see also [8]. As in [1, \S5], let
$R_{FG(n)}^{\rm poly}$ denote the subring of the Green ring $R_{FG(n)}$
spanned by all the isomorphism classes of finite-dimensional
polynomial modules. It is clear from (3.12) that,
for each $r$, the Lie resolvent $\Phi^r$ restricts to a map
$\Phi^r : R_{FG(n)}^{\rm poly}
\to R_{FG(n)}^{\rm poly}$.

Let $t_1,\dots,t_n$ be indeterminates and
let $\Delta$ be the subring of ${\mathbb Z}[t_1,\dots,t_n]$ consisting
of all symmetric polynomials. For each positive integer $r$, the
endomorphism of ${\mathbb Z}[t_1,\dots,t_n]$ satisfying $t_i \mapsto t_i^r$
for all $i$ restricts to an endomorphism of $\Delta$, which we
denote by $\chi^r$.
As explained in [1, \S5], there is a ring
homomorphism ${\rm ch} : R_{FG(n)}^{\rm poly} \to \Delta$ such that,
for every finite-dimensional polynomial $FG(n)$-module $U$,
${\rm ch}(U)$ is the formal character of $U$.

As stated in [6, \S3.2], if $V$ is the natural $FG(n)$-module then
\begin{equation}
{\rm ch}(L^r(V)) = \frac{1}{r} \sum_{d \mid r} \mu(d)
(t_1^d + \cdots + t_n^d)^{r/d}.
\end{equation}
(The right-hand side is formally an element of ${\mathbb Q}
\otimes_{\mathbb Z} \Delta$, but is found to belong to $\Delta$.)
Suppose that $U$ is a finite-dimensional polynomial $FG(n)$-module.
Then we may write ${\rm ch}(U) = w_1 + \cdots + w_m$ where
$m = \dim U$ and $w_1,\dots,w_m$ are monomials in
$t_1,\dots,t_n$. We may choose a basis of $U$ consisting of
elements from weight spaces.  Then every diagonal element of
$G(n)$ is represented on $U$ by a diagonal element of ${\rm GL}(m,F)$.
By (4.1) with $m$ instead of $n$, we obtain
$${\rm ch}(L^r(U)) = \frac{1}{r} \sum_{d \mid r} \mu(d)
(w_1^d + \cdots + w_m^d)^{r/d}.$$
Hence
\begin{equation}
{\rm ch}(L^r(U)) = \frac{1}{r} \sum_{d \mid r} \mu(d)
\chi^d({\rm ch}(U^{r/d})).
\end{equation}

\vspace{0.2in}
\noindent
{\bf Lemma 4.1} \hspace{0.1in}
{\it Let\/ $U$ be a finite-dimensional polynomial
$FG(n)$-module.}\\
(i) {\it For all\/ $r \geq 1$, ${\rm ch}(\Phi^r(U)) = \mu(r)
\chi^r({\rm ch}(U))$.}\\ (ii) {\it If\/ $(r,s) = 1$ then
$${\rm ch}(\Phi^{rs}(U)) =
{\rm ch}(\Phi^r(\Phi^s(U))).$$}

\vspace{0.01in}
\noindent
{\bf Proof:} \hspace{0.1in}  By (3.11),
$${\rm ch}(L^r(U)) = \frac{1}{r} \sum_{d \mid r} {\rm ch}(\Phi^d
(U^{r/d})).$$
Comparing this with (4.2) and using induction on $r$ gives (i).
Hence, for all $W \in R_{FG(n)}^{\rm poly}$,
\begin{equation}
{\rm ch}(\Phi^r(W)) = \mu(r) \chi^r ({\rm ch}(W)).
\end{equation}
Therefore
\begin{eqnarray*}
{\rm ch}(\Phi^r(\Phi^s(U))) & = & \mu(r)
\chi^r ({\rm ch}(\Phi^s(U))) \\
& = & \mu(r) \chi^r (\mu(s) \chi^s ({\rm ch}(U))).
\end{eqnarray*}
However, $\chi^r \circ \chi^s = \chi^{rs}$ and,
if $r$ and $s$ are coprime, $\mu(r)\mu(s) = \mu(rs)$. Therefore
$${\rm ch}(\Phi^r(\Phi^s(U))) =
\mu(rs) \chi^{rs} ({\rm ch}(U)).$$
Hence, by (4.3), ${\rm ch}(\Phi^r(\Phi^s(U))) =
{\rm ch}(\Phi^{rs}(U))$, as required for (ii).
\hfill{$\square$}

\vspace{0.2in}
We recall the Decomposition Theorem (Theorem 2.1). 
The submodules $B_r$ are not 
uniquely determined by $V$, as easy examples show. However, 
it follows from equations (2.1),
by induction, that the $B_r$ are uniquely
determined up to isomorphism. Thus they give uniquely determined
elements of $R_{FG}$.
We shall deduce the following result from the Decomposition Theorem.

\vspace{0.2in}
\noindent
{\bf Theorem 4.2} (Factorisation Theorem) \hspace{0.1in}
{\it Let\/ $F$ be a field of prime
characteristic\/ $p$. Let\/ $G$ be a group. For every non-negative
integer\/ $m$ and every positive integer\/ $k$ not divisible by\/ $p$,
\begin{equation}
\Phi^{p^mk} = \Phi^{p^m} \circ \Phi^k
\end{equation}
and, for every finite-dimensional\/ $FG$-module $V$, the
elements\/ $B_r$ of\/ $R_{FG}$ obtained from the Decomposition Theorem
satisfy
\begin{equation}
p^mB_{p^mk} + p^{m-1}B_{p^{m-1}k}^p + \cdots +
pB_{pk}^{p^{m-1}}
+ B_k^{p^m} = L^k(V^{p^m}).
\end{equation}
}

\vspace{0.01in}
\noindent
{\bf Proof:} \hspace{0.1in}
When we wish to show the dependence on $V$ we write $B_r$ as $B^r(V)$,
and when we wish to show the role of $F$ and $G$ in $\Phi^r$ we
write it as $\Phi_{FG}^r$. 
We first reduce to the case where $F$ is infinite. Let $E$ be an infinite
extension field of $F$. Assume that the theorem holds with
$E$ in place of $F$. Each $FG$-module $V$ determines an $EG$-module
$E \otimes_F V$ by extension of scalars. Thus we obtain a ring
homomorphism $\iota : R_{FG} \to R_{EG}$. It follows from the
Noether--Deuring Theorem (see [5, (29.7)]) that
$\iota$ is an embedding. By [3, Lemma 2.4],
for all $r \geq 1$ and all $V \in {\rm Mod}(FG)$, we have 
$L^r(\iota(V)) = \iota(L^r(V))$ and
$\Phi_{EG}^r \circ \iota = \iota \circ \Phi_{FG}^r$.

By (4.4) over $E$, $\Phi_{EG}^{p^mk} \circ \iota =
\Phi_{EG}^{p^m} \circ \Phi_{EG}^k \circ \iota$ and hence
$\iota \circ \Phi_{FG}^{p^mk} = \iota \circ \Phi_{FG}^{p^m} \circ
\Phi_{FG}^k$. Since $\iota$ is an embedding we obtain (4.4) over $F$.
Let $V \in {\rm Mod}(FG)$.  Applying $\iota$ to
equations (2.1) over $F$ we find that $\iota(B^r(V)) =
B^r(\iota(V))$ for all $r$. By (4.5) over $E$ for the module $\iota(V)$,
$$p^mB^{p^mk}(\iota(V)) + p^{m-1}(B^{p^{m-1}k}(\iota(V)))^p + \cdots
+ (B^k(\iota(V)))^{p^m} = L^k(\iota(V)^{p^m}).$$
Thus
$$\iota(p^mB^{p^mk}(V) + \cdots + (B^k(V))^{p^m}) =
\iota(L^k(V^{p^m})).$$
Since $\iota$ is an embedding we obtain equation (4.5) over $F$.

Hence it is enough to prove the theorem over $E$. In other words we
may assume that $F$ is infinite.

We use induction on $m$ and $k$. If $m=0$ the result is clear because
$\Phi^1$ is the identity map and, for all $V$,
$B^k(V) = L^k(V)$, by (2.1). Hence we may assume that $m \geq 1$ and
that the result holds for $p^id$ with $d$ not divisible by $p$ if $i < m$
or $i = m$ and $d < k$.

Since the functions $\Phi^r$ are linear, in order to prove (4.4)
it suffices to prove
\begin{equation}
\Phi^{p^mk}(V) = (\Phi^{p^m} \circ \Phi^k)(V)
\end{equation}
for all $V \in {\rm Mod}(FG)$. 

Let $V \in {\rm Mod}(FG)$ and 
$n = \dim V$. By choice of a basis of $V$, the representation
of $G$ on $V$ gives a homomorphism $\theta: G \to G(n)$. Furthermore,
$\theta$ induces a ring homomorphism $\theta^* : R_{FG(n)} \to
R_{FG}$. We can regard $V$ as the natural module for $G(n)$, in which
case we write it as $V(n)$. Thus $\theta^*(V(n)) = V$.
We now show that
it is enough to prove the required results for $G(n)$ and $V(n)$.
Suppose that (4.6) and (4.5) hold for $G(n)$ and $V(n)$. By
[3, Lemmas 2.2 and 2.3],
for all $r \geq 1$ and all $U \in {\rm Mod}(FG(n))$, we have
$L^r(\theta^*(U)) = \theta^*(L^r(U))$ and $\theta^* \circ \Phi_{FG(n)}^r
= \Phi_{FG}^r \circ \theta^*$.

Thus, by applying $\theta^*$ to (4.6) for $G(n)$ and $V(n)$, we get
(4.6) for
$G$ and $V$. By applying $\theta^*$ to equations (2.1) for
$G(n)$ and $V(n)$
we get $\theta^*(B^r(V(n))) = B^r(V)$ for all $r$. Hence, by applying
$\theta^*$ to (4.5) for $G(n)$ and $V(n)$, we get (4.5) for $G$ and $V$.

Thus it remains to prove (4.6) and (4.5) for $G(n)$ and $V(n)$.
To simplify
the notation, write $V$ for $V(n)$ 
and $B_r$ for $B^r(V(n))$.

As in \S3, let $\Pi = \Gamma_{FG(n)}[[t]]$. 
Let ${\rm Log}: t\Pi \to t\Pi$ be defined by
$${\rm Log}(f) = - \log(1-f) = f + f^2/2 + f^3/3 + \cdots$$
for all $f \in t\Pi$. Thus ${\rm Log}$ is the inverse of the
function ${\rm Exp}$ defined in \S3.

For $f \in t\Pi$ and $r \geq 1$, we write $f_{(r)}$
for the coefficient of $t^r$ in $f$; thus $f_{(r)} \in \Gamma_{FG(n)}$ 
and $f = \sum_{r \geq 1}f_{(r)} t^r$. 

Let $Q' = \sum_{r \geq 0} {\rm Log}(B_{p^rk}t^{p^r})$,
$Q'' = \sum_{i \geq 0} p^{-i}L^k(V^{p^i})t^{p^i}$ and
$Q = Q' - Q''$. Then it is easily calculated that
$$Q_{(p^i)} = \frac{1}{p^i}B_k^{p^i} + \frac{1}{p^{i-1}}B_{pk}^{p^{i-1}}
+ \cdots + B_{p^ik} - \frac{1}{p^i}L^k(V^{p^i}).$$
By the inductive hypothesis, $Q_{(p^i)} = 0$ for $i < m$. Hence,
by Lemma 3.5, we obtain
\begin{eqnarray*}
{\bf L}^+(Q)_{(p^m)} & = & \Phi^1(Q_{(p^m)}) \; = \; Q_{(p^m)} \\
& = & \frac{1}{p^m}B_k^{p^m} + \frac{1}{p^{m-1}}B_{pk}^{p^{m-1}}
+ \cdots + B_{p^mk} - \frac{1}{p^m}L^k(V^{p^m}).
\end{eqnarray*}
Since ${\bf L}^+$ is additive and 
${\bf L}^+ \circ {\rm Log} = {\cal L}$, by (3.9), we have
$${\bf L}^+(Q') = \sum_{r \geq 0} {\cal L}(B_{p^rk}t^{p^r}).$$
Thus, by (3.5),
$${\bf L}^+(Q')_{(p^m)} = L^{p^m}(B_k) + L^{p^{m-1}}(B_{pk}) +
\cdots + L^1(B_{p^mk}).$$
Therefore, by the Decomposition Theorem,
$${\bf L}^+(Q')_{(p^m)} = L^{p^mk}(V).$$
By Lemma 3.5,
$${\bf L}^+(Q'')_{(p^m)} = \frac{1}{p^m} \sum_{d \mid p^m}
\Phi^d(L^k(V^{p^m/d})).$$
Hence, by (3.11),
\begin{eqnarray*}
{\bf L}^+(Q'')_{(p^m)} & = & \frac{1}{p^m} \sum_{d \mid p^m}
\Phi^d \biggl( \frac{1}{k} \sum_{e \mid k} \Phi^e(V^{p^mk/de})\biggr) \\
& = & \frac{1}{p^mk} \sum_{d \mid p^m,\, e \mid k}
(\Phi^d \circ \Phi^e)(V^{p^mk/de}).
\end{eqnarray*}
We write
$$\Phi^{p^m} \circ \Phi^k = \Phi^{p^mk} + (\Phi^{p^m} \circ \Phi^k
- \Phi^{p^mk})$$
and use the inductive hypothesis to write $\Phi^d \circ \Phi^e =
\Phi^{de}$ for $de < p^mk$. Thus
\begin{eqnarray*}
{\bf L}^+(Q'')_{(p^m)} & = & \frac{1}{p^mk} 
\sum_{s \mid p^mk} \Phi^s(V^{p^mk/s})
+ \frac{1}{p^mk}(\Phi^{p^m}(\Phi^k(V)) - \Phi^{p^mk}(V)) \\
& = & L^{p^mk}(V) + \frac{1}{p^mk}(\Phi^{p^m}(\Phi^k(V)) -
\Phi^{p^mk}(V)).
\end{eqnarray*}

Since ${\bf L}^+(Q)_{(p^m)} = {\bf L}^+(Q')_{(p^m)} -
{\bf L}^+(Q'')_{(p^m)}$, we obtain
$$\frac{1}{p^m}B_k^{p^m} + \frac{1}{p^{m-1}} B_{pk}^{p^{m-1}}
+ \cdots +
B_{p^mk} - \frac{1}{p^m}L^k(V^{p^m}) = 
\frac{1}{p^mk}(\Phi^{p^mk}(V) - \Phi^{p^m}(\Phi^k(V))).$$
Therefore
$$\frac{1}{k}(\Phi^{p^mk}(V) - \Phi^{p^m}(\Phi^k(V)))
= W,$$
where
$$W = p^mB_{p^mk} + p^{m-1}B_{p^{m-1}k}^p + \cdots + pB_{pk}^{p^{m-1}}
+ B_k^{p^m} - L^k(V^{p^m}). $$
By Lemma 4.1, ${\rm ch}(\Phi^{p^mk}(V) - \Phi^{p^m}(\Phi^k(V))) = 0$.
Thus ${\rm ch}(W) = 0$.  By the Decomposition Theorem, 
the modules $B_{p^mk}$, $B_{p^{m-1}k}^p$, \dots, $B_k^{p^m}$
are isomorphic to direct summands of the tensor power $V^{p^mk}$. 
Hence they are tilting modules: see [1, \S5].  Since $p \nmid k$,
$L^k(V^{p^m})$ is also a direct summand of $V^{p^mk}$ (see, for example,
[6, \S3.1]). Thus it too is a tilting module. 
However, tilting modules are determined up to
isomorphism by their formal characters (see, for example, [6]).
Since ${\rm ch}(W) = 0$ it follows that $W = 0$. Thus
$\Phi^{p^mk}(V) - \Phi^{p^m}(\Phi^k(V)) = 0$.
This gives (4.5) and (4.6). \hfill{$\square$}

\vspace{0.2in}
Note that equations (4.5) give a recursive description, within
$\Gamma_{FG}$,
of $B_{p^mk}$ as a polynomial in $L^k(V)$, $L^k(V^p)$, \dots,
$L^k(V^{p^m})$. The polynomials in $B_k, B_{pk}, B_{p^2k},\dots$
occurring on the left-hand side of (4.5) may be recognised as
the Witt polynomials associated to the prime $p$: see
[16, Chapter II, \S6] or [11, Chapter 3].
Thus, in the terminology often used in relation to Witt vectors, 
$L^k(V)$, $L^k(V^p)$, $L^k(V^{p^2})$,
\dots\ are the `ghost' components of $(B_k,B_{pk},B_{p^2k},\dots)$.
However, we have not yet been able to deduce from this more explicit
information on the modules $B_{p^mk}$.

Of particular interest is the case
where $V$ is the natural module for ${\rm GL}(n,F)$. Then, since the
modules $B_k,B_{pk},\dots\,$ and $L^k(V),L^k(V^p),\dots\,$ are
direct summands of tensor powers of $V$, they are determined, up
to isomorphism, by 
their formal characters. We may therefore translate equations (4.5)
into equations in symmetric functions. Closely related equations have
been considered by Reutenauer [14] and Scharf and Thibon [15].

Statement (4.4) in the Factorisation Theorem
is analogous to the property of Adams operations given by
the second part of Lemma 3.2. The Factorisation Theorem and
Lemma 3.4 reduce the study of Lie
resolvents to the study of Adams operations and $p$-power Lie
resolvents.

\vspace{0.4in}
\noindent
{\bf 5. Symmetric powers and Lie powers revisited}
\renewcommand{\theequation}{5.\arabic{equation}}
\setcounter{equation}{0}

\noindent
Let $F$ be a field and $G$ a group. For each positive integer $r$
we may define a function $\rho^r : \Gamma_{FG} \to \Gamma_{FG}$ by
\begin{equation}
\rho^r = \frac{1}{r} \sum_{d \mid r} \Phi^d \circ \psi^{r/d}.
\end{equation}
Thus, since $\psi^1$ is the identity map, the Lie resolvents $\Phi^r$
may be expressed recursively
in terms of the functions $\rho^r$ and the Adams operations. Hence, 
if we assume the Adams operations, knowledge of the functions
$\rho^r$ is equivalent to knowledge of the Lie resolvents.

\vspace{0.2in}
\noindent
{\bf Theorem 5.1} \hspace{0.1in}
{\it Let\/ $G$ be a group and\/ $F$ a field of
prime characteristic\/ $p$. Then\/ $\rho^r = 0$ unless\/ $r$ is a
power of\/ $p$.}

\vspace{0.2in}
\noindent
{\bf Proof:} \hspace{0.1in}
Write $r = p^mk$ where $m \geq 0$ and $k$ is not divisible
by $p$. By
(5.1) and the Factorisation Theorem, 
\begin{eqnarray*}
\rho^r & = & \frac{1}{r} \sum_{d \mid r} \Phi^d \circ \psi^{r/d} \\
& = & \frac{1}{r} \sum_{0 \leq i \leq m} \sum_{e \mid k}
\Phi^{p^ie} \circ \psi^{p^{m-i}k/e} \\
& = & \frac{1}{r} \sum_{0 \leq i \leq m} \sum_{e \mid k} \Phi^{p^i}
\circ \Phi^e \circ \psi^{p^{m-i}k/e}.
\end{eqnarray*}
Hence, by Lemmas 3.4 and 3.2,
\begin{eqnarray*}
\rho^r & = & \frac{1}{r} \sum_{0 \leq i \leq m} \sum_{e \mid k}
\mu(e)(\Phi^{p^i} \circ \psi^e \circ \psi^{p^{m-i}k/e}) \\
& = & \frac{1}{r} \sum_{0 \leq i \leq m} \sum_{e \mid k}
\mu(e)(\Phi^{p^i} \circ \psi^{p^{m-i}k}).
\end{eqnarray*}
If $r$ is not a power of $p$ then $k > 1$ and so $\sum_{e \mid k}
\mu(e) = 0$, which gives $\rho^r = 0$. \hfill{$\square$}

\vspace{0.2in}
\noindent
{\bf Lemma 5.2} \hspace{0.1in}
{\it We have\/ ${\bf L}^* \circ {\bf S}^* = {\bf L}^+ \circ {\bf S}^+$ and,
for all\/ $A \in \Gamma_{FG}$,
$$({\bf L}^* \circ {\bf S}^*)(At) = \rho^1(A)t + \rho^2(A)t^2 + \cdots
\: .$$}

\vspace{0.01in}
\noindent
{\bf Proof:} \hspace{0.1in}
This result is essentially the same as [2, Lemma 4.1], but,
for completeness, we give a proof.

By (3.3) and (3.8), ${\bf L}^* \circ {\bf S}^* = {\bf L}^+ \circ
{\bf S}^+$. Also, by (3.4) and Lemma 3.5,
\begin{eqnarray*}
{\bf L}^+({\bf S}^+(At)) & = & {\bf L}^+(\psi^1(A)t + \frac{1}{2}
\psi^2(A)t^2 + \cdots) \\
& = & \sum_{r \geq 1} \biggl( \sum_{d \mid r} \frac{1}{d}
\Phi^d((d/r)\psi^{r/d}(A))\biggr) t^r.
\end{eqnarray*}
Hence, by the linearity of $\Phi^d$ and (5.1), 
$$({\bf L}^* \circ {\bf S}^*)(At) = ({\bf L}^+ \circ {\bf S}^+)(At)
= \sum_{r \geq 1} \rho^r(A)t^r.$$
{}\hfill{$\square$}

\vspace{0.2in}
Recall that $\Pi = \Gamma_{FG}[[t]]$.
An element of $t\Pi$ will be called $p$-{\it typical} if, for every
positive integer $r$, the coefficient of $t^r$ is zero unless $r$
is a power of $p$.

\vspace{0.2in}
\noindent
{\bf Theorem 5.3} \hspace{0.1in}
{\it Let\/ $G$ be a group and\/ $F$ a field of
prime characteristic\/ $p$. Then\/ $({\bf L}^* \circ {\bf S}^*)(f)$
is\/ $p$-typical for every\/ $p$-typical element\/ $f$ of\/ $t\Pi$.}

\vspace{0.2in}
\noindent
{\bf Proof:} \hspace{0.1in}
By Lemma 5.2, ${\bf L}^* \circ {\bf S}^* = {\bf L}^+ \circ
{\bf S}^+$. By Theorem 5.1 and Lemma 5.2,
$({\bf L}^+ \circ {\bf S}^+)(At)$ is $p$-typical for all
$A \in \Gamma_{FG}$. By Lemmas 3.1 and 3.3, ${\bf L}^+ \circ {\bf S}^+$
has the properties of ${\bf S}^+$ stated in Lemma 3.1. It follows that
$({\bf L}^+ \circ {\bf S}^+)(f)$ is $p$-typical for every $p$-typical
element $f$ of $t\Pi$. \hfill{$\square$}

\vspace{0.2in}
The functions $\rho^r$ seem to have some importance in the
study of Lie powers. As we have already seen, they may be used instead
of the Lie resolvents, but, on the available evidence, their properties
seem to be smoother. Theorem 5.1 shows that $\rho^r = 0$ unless $r$ is
a power of $p$. It is remarkable also that, in the cases which have so
far been calculated, even the functions $\rho^1$, $\rho^p$,
$\rho^{p^2}$, \dots\ are well behaved. Of course, $\rho^1$ is the
identity function. It follows from [2, Corollary 4.5 and Lemma 4.6]
that if $G$ is cyclic of order
$p$ (and $F$ has characteristic $p$) then $\rho^{p^m} = 0$ for all
$m > 1$. In fact, from the results in [2] and [3] it is not
difficult to deduce the same fact for every finite group $G$
such that $p^2 \nmid |G|$. It is interesting to speculate on how far this
generalises to other groups: perhaps, if $G$ is finite, we have 
$\rho^{p^m} = 0$ for all sufficiently large $m$. 

We have used the Factorisation Theorem in the proof of Theorem 5.1 to
obtain the fact that $\rho^r = 0$ unless $r$ is a power of $p$.
However, by [2, Lemma 5.1 (ii)], this fact itself yields
the Factorisation Theorem: hence it is `equivalent' to the Factorisation
Theorem.

We conclude with a few remarks about the easier case where $F$ has
characteristic $0$. In this case, by Lemma 3.4, $\Phi^r = \mu(r)\psi^r$
for all $r$. By a simplified version of the calculation of Theorem 5.1
we obtain $\rho^r = 0$
for all $r > 1$. Thus, by Lemma 5.2, $({\bf L}^+ \circ {\bf S}^+)(At)
= At$ for all $A \in \Gamma_{FG}$. By the properties of ${\bf S}^+$
and ${\bf L}^+$ it follows that ${\bf L}^+ \circ {\bf S}^+$ is the
identity function on $t\Pi$. A similar calculation shows that
${\bf S}^+ \circ {\bf L}^+$ is the identity function. However,
${\bf L}^* \circ {\bf S}^* = {\bf L}^+ \circ {\bf S}^+$ and
$${\bf S}^* \circ {\bf L}^* = {\rm exp} \circ {\bf S}^+ \circ
{\bf L}^+ \circ {\rm log}.$$
Thus we have the following result, closely related to results of
Joyal [12] and Reutenauer [14].

\vspace{0.2in}
\noindent
{\bf Theorem 5.4} \hspace{0.1in}
{\it Let\/ $G$ be a group and\/ $F$ a field of
characteristic\/ $0$. Then\/ ${\bf L}^* \circ {\bf S}^*$ is the
identity function on\/ $t\Pi$ and\/ ${\bf S}^* \circ {\bf L}^*$ is
the identity function on\/ $1 + t\Pi$.}

\vspace{0.2in}
Recall that, by (3.7), ${\bf L}^*(f) = -{\cal L}(1-f)$ for all
$f \in 1 + t\Pi$. 
Since ${\bf S}^* \circ {\bf L}^*$ is the identity function in
characteristic $0$, it follows easily that
$${\bf S}^*({\cal L}(g)) = (1-g)^{-1} = 1 + g + g^2 + \cdots$$
for all $g \in t\Pi$. This may be regarded as a version of the
Poincar\'e--Birkhoff--Witt Theorem (by taking $g = Vt$) and it is 
similar to a result of
Joyal [12, Chapter 4, Proposition 1]
(see also [14, Lemma 3.2]). Similarly, the fact that
${\bf L}^* \circ {\bf S}^*$ is the identity function is a version
of a result of Reutenauer [14, Theorem 3.1].

\vspace{0.4in}
\noindent
{\bf References}

\newcounter{euro}
\begin{list}{\arabic{euro}.\hspace{0.1in}}{\parsep=0pt
\itemsep=0pt \usecounter{euro}}

\item {R. M. Bryant}, ``Free Lie algebras and Adams operations",
{\it J.\ London Math. Soc. (2)\/} {\bf 68} (2003), 355--370.
\item {R. M. Bryant}, ``Modular Lie representations of groups of
prime order", {\it Math. Z.} {\bf 246} (2004), 603--617.
\item {R. M. Bryant}, ``Modular Lie representations of finite groups",
{\it J.\ Austral. Math. Soc.} {\bf 77} (2004), 401--423.
\item {R. M. Bryant and M. Schocker}, ``The decomposition of
Lie powers", preprint \texttt{arXiv:math.RT/0505325}. 
\item {C. W. Curtis and I. Reiner}, {\it Representation Theory of
Finite Groups and Associative Algebras\/}, Wiley--Interscience, New York,
1962.
\item {S. Donkin and K. Erdmann}, ``Tilting modules, symmetric functions,
and the module structure of the free Lie algebra", {\it J.\ Algebra\/} 
{\bf 203} (1998), 69--90.
\item {K. Erdmann and M. Schocker}, ``Modular Lie powers and the Solomon
descent algebra", preprint \texttt{arXiv:math.RT/0408211}.
\item {J. A. Green}, {\it Polynomial Representations of\/ ${\rm GL}_n$},
{\it Lecture Notes in Mathematics\/} {\bf 830}, Springer, Berlin, 1980.
\item {S. Guilfoyle}, {\it On Lie Powers of Modules for Cyclic Groups\/} 
Ph.\ D. thesis, Manchester, 2000.
\item {S. Guilfoyle and R. St\"ohr}, ``Invariant bases for free
Lie algebras", {\it J. Algebra\/} {\bf 204} (1998), 337--346.
\item {M. Hazewinkel}, {\it Formal Groups and Applications\/},
Academic Press, New York, 1978.
\item {A. Joyal}, ``Foncteurs analytiques et esp\`eces de structures",
in {\it Combinatoire \'Enum\'erative\/}, edited by G. Labelle and P. Leroux,
{\it Lecture Notes in Mathematics\/} {\bf 1234}, Springer, Berlin, 1986,
pp.\ 126--159.
\item {C. Reutenauer}, {\it Free Lie Algebras\/}, Clarendon Press, Oxford,
1993.
\item {C. Reutenauer}, ``On symmetric functions related to Witt
vectors and the free Lie algebra", {\it Adv.\ Math.} {\bf 110} (1995),
234--246.
\item T. Scharf and J.-Y. Thibon, ``On Witt vectors and symmetric
functions", {\it Algebra Colloq.} {\bf 3} (1996), 231--238. 
\item {J.-P. Serre}, {\it Local Fields\/}, Springer, New York, 1979.
\item {R. St\"ohr}, ``Restricted Lazard elimination and modular Lie
powers", {\it J. Austral.\ Math.\ Soc.} {\bf 71} (2001), 259--277. 

\end{list}

\end{document}